\documentclass{amsart}

\usepackage{amsmath,amssymb,amsthm,epsfig,psfrag}

\newtheorem{theorem}{Theorem}
\newtheorem{proposition}[theorem]{Proposition}
\newtheorem{lemma}[theorem]{Lemma}

\theoremstyle{definition}
\newtheorem{definition}[theorem]{Definition}

\def\cal{\mathcal}

\def\Bbb{\mathbb}
\def\bar{\overline}

\def\Z{\Bbb{Z}}
\def\N{\Bbb{N}}

\def\Reals{\Bbb{R}}

\def\ni{\noindent}

\def\ssm{\smallsetminus}

\def\ms{\medskip}

\def\onto{{\kern3pt\to\kern-8pt\to\kern3pt}}
\def\<{\langle}
\def\>{\rangle}
\def\|{{\ |\ }}

\def\0{\mathbf{0}}

\def\II{\cal I}

\def\AA{\cal A}
\def\BB{\cal B}
\def\CC{\cal C}

\def\SS{\cal S}
\def\DD{\cal D}
\def\LL{\cal L}
\def\PP{\cal P}

\newcommand{\set}[1]{\left\{#1\right\}}

\newcommand{\abs}[1]{\left|#1\right|}

\renewcommand{\ni}{\noindent}
\renewcommand{\ss}{\smallskip}
\renewcommand{\ms}{\medskip}
\newcommand{\bs}{\bigskip}
\def\*{^{\star}}

\begin{document}

\markboth{Tim R.\ Riley and Andrew D.\ Warshall}{The unbounded dead-end depth property is not a group invariant}

\title{THE UNBOUNDED DEAD-END DEPTH PROPERTY \\ IS NOT A GROUP INVARIANT}

\author{TIM R.\ RILEY}

\author{ANDREW D.\ WARSHALL}

\date{31 March 2005, revised 9 September 2005}

\begin{abstract}
The \emph{dead-end depth} of an element $g$ of a group with finite
generating  set $\AA$ is the distance from $g$ to the complement of the radius
 $d_{\AA}(1,g)$ closed ball, in the word metric $d_{\AA}$.  We exhibit a
 finitely presented group $K$ with two finite generating sets $\AA$ and
$\BB$  such that dead-end depth is unbounded on $K$ with respect to $\AA$ but is
bounded above by two with respect to $\BB$.

\ms

\ni Keywords: dead-end depth, lamplighter.
\end{abstract}

\maketitle

\section{Introduction} \label{intro} 

Suppose $G$ is a group with finite generating set $\AA$ and associated word metric
 $d_{\AA}$.
The {\em dead-end depth} (or, more concisely, the
\emph{depth}) of $g \in G$ with respect to $\AA$ is the distance between $g$ and the complement of the ball of radius $d_{\AA}(1,g)$ centered at $1$ in $G$.  (If the ball is all of $G$ then
the depth of $g$ is infinite.)  An element $g$ of depth greater than one is
called a \emph{dead end} because a geodesic from $1$ to $g$ in the Cayley graph
cannot be continued to a group element beyond $g$. 

In general, depth depends on the generating set.  
For example, all elements of $\Z = \langle a \rangle$ have depth $1$ in
$(\Z,\set{a})$, but the depth of $a$ in $(\Z, \set{a^2, a^3})$ is  $2$.  This article 
addresses how radical the dependence can be.  

Cleary \& Taback \cite{CT1,CT2} showed depth to be unbounded in the
lamplighter group $\Z_2\wr\Z$ and similar wreath products, with
respect to certain finite generating sets -- for instance $\set{a,t}$, where
$\Z_2\wr\Z$ is presented by $\langle a,t\mid
a^2;[a^{t^i},a],\forall i\in\Z\rangle$. (Independently, Erschler observed that
$\Z_2\wr\Z$  provides an example resolving the closely
related Question~8.4 of Bowditch in \cite{Bestvina}.) This prompted
the question (asked by Taback, lecturing at CUNY in April 2004, and
subsequently in print by Cleary and the first author in \cite{CR}) whether the
property of depth being unbounded is a \emph{group invariant}, that
is, is independent of the finite generating set.  Indeed, there was
speculation on whether the property might be a quasi-isometry
invariant.  We answer these questions negatively. (We denote the commutator $a^{-1} b^{-1} a b$ by
 $[a,b]$ and the conjugates $b^{-1}ab$ and $b^{-1}a^{-1}b$ by $a^b$ and $a^{-b}$, respectively.)

\begin{theorem} \label{Thm K} The group $K$ with
finite presentation
$$\langle \ a,s,t,u \ \mid \ 
a^2, \, [a,a^t], \, [s,t], \, [s,u], \, [t,u], \,
 {a^{-s}}aa^t, \, {a^{-u}}a^t \
 \rangle$$
has unbounded dead-end depth with respect to the generating set
\[
\AA \ = \ \{ \,
a,s,t,u,as,at,au,sa,ta,ua,asa,ata,aua \, \}
\]
but dead-end depth bounded above by $2$
with respect to 
\begin{eqnarray*}
\BB \ = \
 \AA  & \cup & 
\{ \, tu,atu,tau,tua,atau,atua,taua,ataua, \\ & &
 ut^{-1},aut^{-1},tat^{-1},t^{-1}at,atat^
{-1},at^{-1}at \, \}.
\end{eqnarray*}
\end{theorem}

\bs

Dead ends and depth have been studied in a variety of settings.  In 
addition to those already mentioned these include $\textup{SL}_2(\Z)$ and $\langle
x,y \mid x^3, y^3, (xy)^k \rangle$ (see \cite{Bogop}), Thompson's group $F$ (see
 \cite{Fordham}), and presentations satisfying the $C'(1/6)$ small
cancellation condition (see \cite{Champ}).  Also IV.A.13,14 of  \cite{dlH} contain a
 discussion of (non-) dead ends under the name ``extension property 
for  geodesic segments.''

Dead ends and depth (not so-termed in \cite{LPP}) are the key to 
a striking result of Lyons, Pemantle \& Peres \cite{LPP}: they show that random walks on the lamplighter group biased suitably towards the identity (``homesick'' random walks) move outward from the identity faster than simple random walks.  The point is that the homesick random walk
will escape dead ends faster than a simple random walk. 

In general, it seems hard to
understand the behaviour of depth in a group as the
 finite generating set varies.  Even depth in $\Z$ is not entirely
 straightforward.  For any given finite generating set, depth is bounded -- indeed, Bogopol'ski{\u\i} showed this to be the case to be the case in all infinite hyperbolic groups \cite{Bogop}. However, as we will show in Section~\ref{depth in Z}, depth
 in $\Z$ is not uniformly bounded as the generating set varies:

\begin{proposition} \label{depth in Z result} For all $k \in \N$, there
exists a finite generating set $\AA$ for $\Z$ for which there is a group element of depth greater than $k$ in $(\Z,\AA)$.
\end{proposition}

\bs
En route to Theorem~\ref{Thm K} we will prove an analogous result for a group that is finitely generated but not finitely presentable:
\ms

\begin{theorem}  \label{Thm H}
The group $H$ presented by
$$ \langle \ a, t, u  \mid a^2, \ [t,u], \ a^{-u}a^t\,; \ \forall i \in 
\Z, \
 [a,a^{t^i}] \ \rangle$$
has unbounded depth with respect to the generating set
$$ \CC =  \set{ \, a, t, u, at, ta, ata, au, ua, aua \, }$$ but depth
 bounded above by $2$ with respect to
\begin{eqnarray*} \DD \ = \ \CC  & \cup &  \{ \, tu, atu, tau,
tua, atau, atua, taua, ataua, \\ & & \ ut^{-1}, aut^{-1}, tat^{-1}, t^{-1}at,
 atat^{-1}, at^{-1}at \, \}. \end{eqnarray*} 
 \end{theorem}

\bs
The construction of $H$ is similar to that of $\Z_2\wr\Z$:
 specifically, $H = \left( \bigoplus_{i \in \Z} \Z_2 \right) \rtimes
\Z^2,$ where  $a_i$
 generates the $i$-th copy of $\Z_2$ in the direct sum, $\Z^2 =
 \langle t,u
 \rangle$, and the actions of $u$ and $t$ are given by ${a_i}^t=
 {a_i}^u=
 a_{i+1}$.  Defining $a:=a_0$ we find $a_{i} = a^{t^i}$, and the
presentation
 given in Theorem~\ref{Thm H} can be obtained by simplifying $$\langle \
(a_i)_{i \in \Z}, t, u \mid \ [t,u] \,; \ \forall i,j \in \Z, \
 [a_i,a_j], \ {a_i}^t{a_{i+1}}^{-1}, \ {a_i}^u{a_{i+1}}^{-1}, \ {a_i}^2 \
 \rangle.$$
 (Identifying  $t$ and $u$ retracts $H$ onto $\Z_2\wr\Z$ and so gives a way of obtaining a presentation of $\Z_2\wr\Z$ from a presentation for $H$ that shows $H$ not to be finitely presentable because $\Z_2\wr\Z$ is not  finitely presentable.)    
 The group $G$ of \cite{CR} presented by
$$\langle \ a,s,t \ \mid \ a^2, \ [a,a^t], \ [s,t], \ a^{-s}aa^t \
\rangle$$
was the first example of a \emph{finitely presentable} group 
with unbounded
 dead-end depth with respect to some finite generating set, specifically 
$\set{a,s,t,at,ta,ata,as,sa,asa}$. The group $K$ of Theorem~\ref{Thm 
K} is related to $G$ similarly to how $H$ is related to the lamplighter 
group.

\bs

This article is organised as follows. Sections~\ref{lamplighter model
for H} and \ref{lamplighter model
for K} describe \emph{lamplighter models} that aid understanding of the 
geometries of $H$ and $K$, and Sections~\ref{proof of Thm H} and \ref{proof of Thm K}
contain proofs of Theorems~\ref{Thm H} and \ref{Thm K}, respectively.  The brief Section~\ref{depth in Z} contains a proof of Proposition~\ref{depth in Z result} on depth in $\Z$.

\ms

\ni \emph{Acknowledgement.} We thank Joshua Zelinsky for his assistance
with Proposition~\ref{depth in Z result} and an anonymous referee for a careful reading.

\section{A lamplighter model for $H$} \label{lamplighter model for H}

Cannon's lamplighter model is a well-known aid to understanding $\Z_2 \wr \Z$.
 It involves a lamplighter moving along a $\Z$-indexed string of lamps and toggling them between on and off.  We will describe a similar model for $H$.

\begin{figure}[ht] \psfrag{a}{$a$}
\psfrag{t}{$t$} 
\psfrag{u}{$u$}
\centerline{\epsfig{file=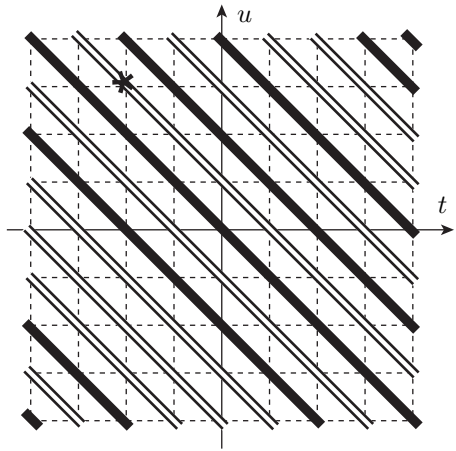}}
\caption{Diagonal strip-lights in the lamp-lighting model for $H$.  Of the
lamps shown, numbers $-7, -5, -4, -3, -1, 1, 3, 5, 6$ (those with white stripes) are illuminated, and the
 lamplighter is at $(-2,3)$.} \label{lamp-grid} \end{figure}

We will give a faithful transitive left action of $H$ on
$\PP_{\!\textit{fin}}(\Z) \times \Z^2$, where
$\PP_{\!\textit{fin}}(S)$ denotes the collection of finite
subsets of a set $S$.  Killing $a$ defines a retraction $\LL:H \onto \langle
t,u \rangle \cong \Z^2$ that gives the location $\LL(g)=(q,r)$ of a
lamplighter among the $\Z^2$-lattice points of the plane after the left
action of $g$ on $(\emptyset, \mathbf{0}\!\,)$. (Here, $q$ and $r$ are
the $t$- and $u$-coordinates, respectively, of the lamplighter and
$\mathbf{0}:=(0,0)$.)  In contrast to the lamplighter model for $\Z_2
\wr \Z$, we conceive the lamps to be infinitely long bulbs (like a long fluorescent tube).  These are $\Z$-indexed and the $i$-th lamp runs along the entire length of the line $x+y=i$ as illustrated in
Figure~\ref{lamp-grid}.  

Defining actions of $t$, $u$ and $a$ on $\PP_{\!\textit{fin}}(\Z) \times \Z^2$ as follows extends to an action in which an element of $H$  represented by a word $w$ on $\set{t,u,a}$ acts by the composition of the actions of the letters of $w$, beginning with the right-most.  (The reader can check that the defining relations act trivially and so this action is well-defined.) 

The actions of $t$ and $u$ are to move the lamplighter one
unit in the $t$- and $u$-directions, respectively. When $(q,r)$ is the
location of the lamplighter the action of $a$ is to toggle the
$(q+r)$-st lamp between on and off.  In other words, $a$ toggles the lamp running through the lamplighter's position.   We define $\II:H \onto \PP_{\!\textit{fin}}(\Z)$ by, given $g\in H$, setting $\II(g)$ to be the set of lamps illuminated after $g$ has acted on $(\emptyset,
\0\!\,)$, and we see that $\LL(g)$ is then the location of the lamplighter.   

Checking that the action is transitive is straightforward and is left to the reader. To
show the action is faithful we suppose $g \in H$ satisfies
$g(\emptyset, \0 \!\,) = (\emptyset, \0\!\,)$ and we check that $g=1$
in $H$. Let $w$ be a word representing $g$.  Reading $w$ from right to
left determines a path followed by the lamplighter in the grid, starting and finishing at $\0$
and in the course of which lamps are toggled between on and off.  The
relation $a^t=a^u$ implies that $a=a^{ut^{-1}}$ and
$a=a^{tu^{-1}}$. These relations can be used to alter $w$ to another
word $w'$ that also represents $g$ in $H$ and that determines a path
in the grid in the course of which lamps are only toggled on and off
when the lamplighter is on the $t$-axis. (For example, if $w$ is
$uatu^{-1}at^{-1}$ then  $w'$ can be taken to be
$uu^{-1}tat^{-1}utu^{-1}at^{-1}$.)  The relation $[t,u]=1$ together
with free reduction can then be used to alter $w'$ to another word
$w''$ that represents $g$ and that determines a path in which the
lamplighter does not leave the $t$-axis. (In our example, this would
give us $w''=ta^2t^{-1}$.)  Identifying $t$ and $u$ defines a retraction
of $H$ onto the lamplighter group, under which $g$ is mapped to an
element also represented by $w''$.  So $g=1$ in $H$ because the action
of the lamplighter group  $\Z_2 \wr \Z$ in Cannon's model is faithful.

\section{Proof of Theorem~\ref{Thm H}} \label{proof of Thm H}

The following proposition establishes that depth is unbounded on $(H,
 \CC)$.  For, suppose we take $g_n \in H$ with $\II(g_n)=\set{-n,n}$ and
$\LL(g_n) = \0$.  %(For example, the word $t^nat^{-2n}at^n$ represents $g_n$.)
Then $d_{\CC}(1,g_n) =4n$ and, with respect to $\CC$, the distance from $g_n$
to the complement of the radius $4n$ closed ball about $1$ is at least $n+1$
since the lamplighter has either to toggle a light outside $\set{-n, \ldots, n}$ or to end up outside $\set{-n,\ldots,n}$ and so in either case
must travel outside $D_n$.

\begin{proposition} Define $$D_n:=\set{(q,r) \in
\Z^2 \mid \abs{q} + \abs{r} \leq n}.$$ If $g \in H$
 satisfies $\LL(g) \in D_n$ and $\II(g) \subseteq \set{-n, \ldots, n}$
then $d_{\CC}(1,g) \leq 4n$.  If, in addition, $\LL(g) = \0$ and $\set{ -n ,n}  \subseteq
 \II(g)$, then $d_{\CC}(1,g) = 4n$. 
 \end{proposition}

\begin{proof}
The salient feature of $\CC$ is that, when the lamplighter moves in the
 grid from a vertex $v_1$ to an adjacent vertex $v_2$, one or both of
 the lights at $v_1$ and $v_2$ can be toggled between on and off with
 no additional cost to word length.  So, for $g \in H \ssm \set{a}$,
 we find $d_{\CC}(1,g)$ is the length of the shortest path in the grid
 that starts at $(0,0)$, visits all of the lights that have to be
 illuminated, and finishes at $\LL(g)$.

Assume $g \in H \ssm\set{a}$ with $(q,r):=\LL(g) \in D_n$ and $\II(g)
\subseteq \set{-n, \ldots, n}$.  Define $$v_1 \ := \
\left(\frac{n+q-r}{2},\frac{n-q+r}{2}\right), \ \ \ v_2 \ := \
\left(\frac{-n+q-r}{2},\frac{-n-q+r}{2}\right),$$ the points of
intersection of the line $x-y=q-r$ with the lines $x+y=n$ and
$x+y=-n$, respectively.  Either both $v_1$ and $v_2$ are in $\Z^2$, or
both are in $(1/2+\Z)^2$.  In the latter case redefine $v_1$ and $v_2$
by adding $(1/2,-1/2)$ to both.  Then (in either case) there is a grid
path from $v_1$ to $v_2$ of length $2n$ that passes through
$(q,r)$.  Assume $(q,r)$ is closer to $v_1$ than $v_2$ along this
path; otherwise interchanging $v_1$ and $v_2$ in the following gives
the required result. There are grid paths of length $n$ from $\0$ to
$v_2$, of length $2n$ from $v_2$ to $v_1$, and of length at most $n$
from $v_1$ to $(q,r)$.  Concatenating gives a path of length at most
$4n$ in the course of which every light in $\set{-n, \ldots, n}$ is
visited.  So $d_{\CC}(1,g) \leq 4n$. If, in addition, $\LL(g) = \0$
and $\set{-n,n} \subseteq \II(g)$ then the lamplighter must visit lamp
$-n$ and then lamp $n$, or vice versa, and then return to $\0$.  Grid
paths of length $n$, $2n$ and $n$, respectively, are both necessary
and sufficient for these three components of the journey.  So
$d_{\CC}(1,g) = 4n$ as required.
\end{proof}

\bs

We will show that depth is identically $1$ in $(H, \DD)$ except at
$a$. The depth at $a$ is $2$ because all group elements at
distance $1$ from it are either other elements of $\DD$ or the identity, while these
other elements of $\DD$ themselves have depth $1$, as we will show.  This will
prove Theorem~\ref{Thm H}.  We begin by defining a $\DD$-path in
$\Reals^2$ to be a concatenation of straight-line \emph{segments} each of
which connects a point in $\Z^2$ to a point in $(1/2+\Z)^2$, has
length $\sqrt{2} /2$, and (so) has slope $\pm 1$.  Define the length
of a $\DD$-path to be the number of such constituent line
segments. The following lemma is straightforward.

\begin{lemma}\label{norm}
For $(q,r) \in \Z^2$ the length of the shortest $\DD$-path from $\0$ to
 $(q,r)$ is $\abs{q+r} + \abs{q-r}$, that is the $\ell_1$-norm of $(q,r)$ 
with respect to the basis
 $\set{(1/2,1/2),(-1/2,1/2)}$.
\end{lemma}

Our next lemma relates $\DD$-paths and the word metric $d_{\DD}$.

\begin{lemma}
For $g \in H \ssm \set{a}$, the distance $d_{\DD}(1,g)$ is half the length 
of the  shortest $\DD$-path from $(0,0)$ to $\LL(g)$ that visits all the lamps in $\II(g)$.
\end{lemma}

\begin{figure}[ht]
\psfrag{tu}{(3) \ $tu, atu, tau, tua, atau, atua, taua, ataua$}
\psfrag{t}{(1) \ $t, ta, at, ata$}
\psfrag{u}{(2) \ $u, ua, au, aua$}
\psfrag{or}{,}
\psfrag{ut^{-1}}{(4) \ $ut^{-1}, aut^{-1}$}
\psfrag{tt2}{(6) \ $tat^{-1}, atat^{-1}$}
\psfrag{tt1}{(5) \ $t^{-1}at, at^{-1}at$}
\psfrag{7}{(7)}
\psfrag{1}{$1$}
\centerline{\epsfig{file=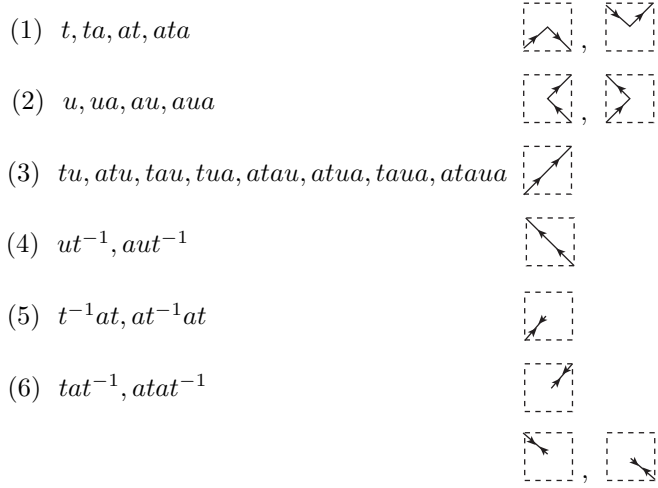}}
\caption{Generators in $\DD$ and the corresponding $\DD$-paths.}
 \label{BB-moves}
\end{figure}

\begin{proof}
Let $\mu$ be a minimal length $\DD$-path from $\0$ to $(q,r)$ 
that, en route, toggles some of the lamps it visits in such a way as to 
illuminate the pattern $\II(g)$.  Choose some of the vertices at the start or end of
segments in $\mu$ to be \emph{distinguished}, with the selection being 
made in such a way that toggling all the lamps incident with distinguished 
vertices lights  the pattern $\II(g)$.

The lamplighter in the grid model for $H$ cannot follow $\mu$ because of being  constrained to move between points in $\Z^2$.  However, $\mu$ has even length as  the points at the ends of segments are alternately in $\Z^2$ and $(1/2+\Z)^2$.

Consider $\mu$ two segments at a time. 
A pair of adjacent segments combine to give a path between two $\Z^2$ 
points (that may include some of the distinguished vertices). Well, a generator in $\DD \ssm \set{a}$  acts to  move the lamplighter between two $\Z^2$ points and possibly toggle some lamps between the initial and the final location of the lamplighter. The construction of $\DD$ is such that elements of $\DD^{\pm 1} \ssm \set {a^{\pm1}}$ are in correspondence with certain pairs of  adjacent length-$\sqrt{2} /2$ diagonal line-segments along which some vertices are
distinguished.  There are $16$ ways to construct a path of two joined diagonal segments, $10$ of which are depicted in Figure~\ref{BB-moves};   the remaining $6$ are the inverses of the paths in the first four rows. Corresponding elements of $\DD$ are shown in the figure, and the locations of the letters $a$ dictate where the incident lamp is to be toggled.  The two cases in row seven cannot occur in $\mu$ because $\mu$ is of minimal length.  For the same reason, the cases in
rows five and six can only occur when the lamp incident with the midpoint of the square is to be toggled.  Only two words occur in the fourth row, one with an $a$ and one without, as the corresponding $\DD$-path is incident with only one lamp.

So we can produce a word $w$ on the alphabet $\DD^{\pm 1}$ such that the length  of $w$ is half that of $\mu$, that, read right to left, describes a grid path from $\0$ to $(q,r)$ and along which lamps are toggled to illuminate the configuration $\II(g)$.  In the same manner, given a geodesic word $w$ on 
$\DD^{\pm 1}$ representing some $g \in H \ssm \set{a}$, we can produce a $\DD$-path of length twice that of $w$ from $\0$ to 
$(q,r)$ that, en route, switches on the lamps $\II(g)$.  (There will be no generators $a$ in $w$ as they could be absorbed into an adjacent generator, reducing word length.)  
\end{proof}

\ss

The two lemmas combine to give:

\begin{proposition} Suppose $g \in H
\ssm \set{a}$.  Let $(q,r):=\LL(g)$.  Then $d_{\DD}(1,g)$ is  half the
 length of the shortest $\DD$-path that starts at $(0,0)$, travels along
the line $x=y$ visiting all lamps in $\II(g)$, then arrives at
 $((q+r)/2,(q+r)/2)$ \emph{(}the point of intersection of the lines
 $x=y$ and $x+y=q+r$\emph{)} and finally travels along $x+y=q+r$ to $(q,r)$.
\end{proposition}

So, given $g \in H \ssm\set{a}$, either $d_{\DD}(1,(ut^{-1})g)=
d_{\DD}(1,g) +1$ or $d_{\DD}(1,(tu^{-1})g)= d_{\DD}(1,g) +1$.  This
completes the proof of Theorem~\ref{Thm H}.

\section{A lamplighter model for $K$} \label{lamplighter model for K}

In Section~\ref{lamplighter model for H} we gave a lamplighter model for $H$ that amounted to a  left action on $\PP_{\!\textit{fin}}(\Z)\times\Z^2$.  Similarly, $K$ has a lamplighter model: a left action on $\PP_{\!\textit{fin}}(\Z)\times\Z^3$ that will be given by two maps, $\II:K\onto\PP_{\!\textit{fin}}(\Z)$ and $\LL: K \onto \Z^3$, such that $g(\emptyset,\0) = (\II(g), \LL(g))$.  This action will be faithful and transitive, or, equivalently, $\II \times \LL$ will be bijective.   

\begin{figure}[ht]
\psfrag{a}{$a$}%
\psfrag{s}{$p$}%
\psfrag{t}{$q$}%
\centerline{\epsfig{file=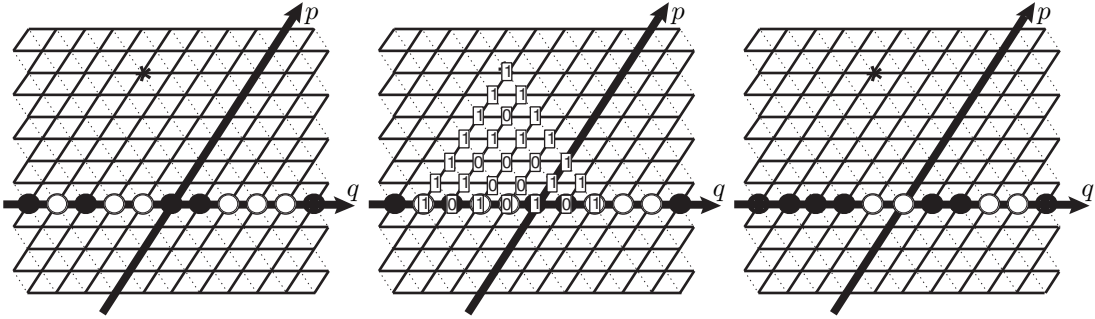}}
 \caption{An example of the action of $a$.  The left diagram shows $g ( \emptyset, \textbf{0} )$ and the right diagram shows $ag ( \emptyset, \textbf{0} )$, where $g=s^6at^{-2}at^{-1}at^{-3}at^{-1}at^{-1}at^4$.  Along the $t$ axis, open circles indicate illuminated lamps and filled-in circles indicate lamps which are off.} \label{lamp-lighting grid}
\end{figure}

Define $\LL : K \onto \langle s,t,u \rangle \cong \Z^3$  to be the retraction given by killing $a$.  This  defines the location $\LL(g)=(p,q,r)$ of the lamplighter among the $\Z^3$-lattice
points of $\Reals^3$ after the action of $g \in K$ on
$(\emptyset,\0 \!\,)$, where $p$, $q$ and $r$ are the $s$-, $t$- and
$u$-coordinates.  The definition of $\II : K \onto \PP_{\!\textit{fin}}(\Z)$ involves the group $G$ presented by 
$$\langle \ a,s,t \ \mid \ a^2, \ [a,a^t], \ [s,t], \ a^{-s}aa^t \ \rangle$$ and related to $K$ as discussed in Section~\ref{intro}. The following lamplighter model for $G$ was given in \cite{CR}.   Envisage a $\Z$-indexed set of lamps that are arranged along the $q$-axis in the $p,q$-plane in which the $p$-axis is skewed so as to make an angle $\pi/3$ with the $q$-axis  --- see Figure~\ref{lamp-lighting grid} (reproduced from \cite{CR}).   Regard an element of $\PP_{\!\textit{fin}}(\Z)\times\Z^2$ as a pair consisting of a finite set of lamps (\emph{illuminated} lamps) together with a lattice point (\emph{lamplighter position}).  Define a (faithful, transitive) left action of $G$ on  $\PP_{\!\textit{fin}}(\Z)\times\Z^2$ by letting $s$ and $t$ move the lamplighter a distance $1$ in the $p$- and $q$-directions, respectively, and by letting $a$ ``press a button'' at the location of the lamplighter that toggles the lights at the locations of the $1$s in a modulo $2$ Pascal's triangle, suspended  from (or growing up from, when below the $q$-axis) the location of the lamplighter, as illustrated in Figure~\ref{lamp-lighting grid}.     
  
\ni  Identifying $t$ and $u$ defines a retraction $\Phi: K \onto G$.  For $g \in K$, define $\II(g)$ to be the $\PP_{\!\textit{fin}}(\Z)$ component of the image of $(\emptyset,\0)$ under the action of $\Phi(g)$.  

The  map $\II \times \LL : K \to \PP_{\!\textit{fin}}(\Z) \times \Z^3$ is easily checked to be surjective.  To show it is injective, suppose $g \in K$ satisfies $(\II \times \LL)(g) =  ( \emptyset, \textbf{0} )$; the action defined above of $G$ on   $\PP_{\!\textit{fin}}(\Z)\times\Z^2$ is faithful \cite{CR} and so $ g  \in \textup{Ker} \, \Phi$.  But $\textup{Ker} \, \Phi$ is the normal closure of $\langle tu^{-1} \rangle$ in $K$. Since $tu^{-1}$ is in the center of $K$, it generates its own normal closure, and so $\textup{Ker} \LL  \cap \textup{Ker} \, \Phi = \set{ 1}$. Hence  $g=1$.  

To understand this action of $K$ geometrically, let $(p,q,r)$ be a co-ordinate system for  $\mathbb{R}^3$ 
and conceive the lamps to be infinitely long bulbs indexed by $\Z$,  the $i$-th bulb running along the line $q+r=i,p=0$ in $\mathbb{R}^3$.  The $q$- and $r$-axes are taken to be mutually orthogonal, but the $p$-axis makes a $\pi/3$ angle with the line $q=r$ and lies in the vertical plane containing it.   
The actions of $s$, $t$ and $u$ are to move the lamplighter one unit in the $p$-, $q$- and $r$-directions, respectively.   Envision all the planes of constant $q-r$ to be tessellated by equilateral triangles in such a way that two of the three sides represent a unit increment of $p$ and $q+r$, respectively, and so that the vertices lie at integer values of $p$ and $q+r$. (Thus the vertices along the line $p=0$ in each plane are at the intersections of the strip-lights with that plane.)   The action of $a$ is to ``press a button'' at $\LL(g)$  toggling the bulbs incident with the locations of the $1$s in a modulo $2$ Pascal's Triangle in a plane of constant $q-r$ suspended (or growing up from if $p<0$) from the location of the lamplighter.

\section{Proof of Theorem~\ref{Thm K}}\label{proof of Thm K}

\begin{figure}[ht]
\psfrag{a}{$a$}
\psfrag{s}{$p$}
\psfrag{t}{$q$}
\centerline{\epsfig{file=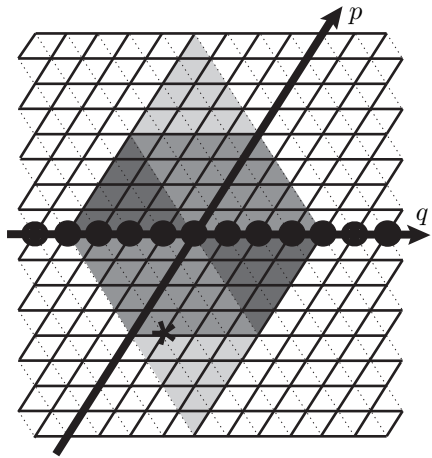}} \caption{Diamond, hexagonal, 
and triangular regions in the $p,q$-plane} \label{grid2}
\end{figure}

We will first show that $(K, \mathcal{A})$ has unbounded depth.  The proof begins with
the following lemma, the analogue of Proposition~2 in \cite{CR}.  Define $H_n$, as pictured in Figure~\ref{grid2} (reproduced from \cite{CR}), to be the subset of $\Z^2$ consisting of lattice points in the (closed) hexagon with corners at $(\pm n, 0)$, $(0, \pm n)$, $(n, -n)$ and $(-n, n)$, and 
 $T_n$ to be the subset of $\Z^2$ consisting of lattice points in the (closed)
 triangle with vertices at $(0,0)$, $(0,-n)$ and $(n,-n)$.  That is, of the shaded regions above the mid-line in the figure, $T_n$ is the darkest. Define
 $$P_n \ := \ \set{ (p,q,r) \mid (p,q+r) \in H_n \textup{ and }
 \abs{q} + \abs{r} \le n}.$$

\begin{lemma}\label{distmax}
All group elements $g$ with $\II(g)\in\set{-n,\dots,n}$ and $\LL(g)\in P_n$ \textup{(}where $n>0$\textup{)} are within $4n$ of the identity in $(K,\AA)$.
\end{lemma}

\begin{proof}
If $\II(g)\subseteq\set{0}$ and $\LL(g)=0$ then $d_{\AA}(1,g) \leq 1 \leq 4n$, as required.

For other $g \in K$ we proceed roughly as in \cite{CR}. As there, the important feature of $\AA$ is that the word length of any $g \in K \ssm \set{a}$ is the length of the shortest path from $\0$ to $\LL(g)$ that visits all the lamps in $\II(g)$. This is because any button encountered en route may be pressed at no additional cost to word length.  

We will only prove the result for $p \geq 0$.  An analogous approach establishes the result for $p <0$.

Suppose $(p,q+r) \in T_n$.  Then, in particular, $p \leq -q -r$.  Let the lamplighter follow the path $\psi$ with straight-line segments:
$$(0,0,0) \ \to \ (0,n,0) \ \to \ (0,q,0) \ \to \  (0,q,-n-q) \ \to \ (0,q,r)  \ \to \  (p,q,r).$$ 
These segments have lengths $n$, $n-q$, $n+q$, $n+q+r$ and $p$, respectively.  (Note that $n-q$, $n+q$ and $n+q+r$ are all non-negative because $\abs{q} + \abs{r} \leq n$.)  So, as $p+q+r \leq 0$,  the total length of $\phi$ is at most $4n$.  

Suppose $(p,q+r) \in H_n \ssm T_n$ and $p\geq 0$.  Then $0 < p+q+r \leq n$.  If $r \geq 0$ then let the path $\psi$ be comprised of the successive straight-line segments
\begin{equation*}
\begin{array}{l}
(0,0,0) \  \to \ (0, - n,0) \ \to \ (0,-p, 0) \ \to \  (p,-p,0) \ \to \ (p,-p,r)  \\  \ \ \ \ \ \  \to \  (p ,n-p-r, r ) \ \to \  (p ,q, r ),\end{array}
\end{equation*}
of lengths $n$, $n-p$, $p$, $r$, $n-r$ and $n-p-q-r$, and so of total length $4n- p-q-r$, which is at most $4n$.  If $r <0$ then take $\psi$ to be
\begin{equation*}
\begin{array}{l}
(0,0,0) \  \to \ (0, - n,0) \ \to \ (0,-p, 0) \ \to \  (p,-p,0) \ \to \ (p,n-p,0)  \\  \ \ \ \ \ \  \to \  (p ,n-p, r ) \ \to \  (p ,q, r ),\end{array}
\end{equation*}
 which has length $n+ (n-p) + p + n-r+ \abs{n-p-q} = 3n-r+ \abs{n-p-q}$.  This is at most $4n$ because if $n-p-q \geq 0$  then it equals $4n-p-q-r$, and if  $n-p-q \leq 0$ then it equals $2n +p+q-r$ and $p \leq n$ and $q-r \leq n$ as $\abs{q}-r = \abs{q}+\abs{r} \leq n$.
 
In every case, adding the second and third co-ordinates projects $\psi$ onto the path $\phi$ of the proof of Proposition~2 in \cite{CR}.  It is shown there that given any set of bulbs with numbers in $\set{-n, \ldots, n}$, there is a combination of buttons on $\phi$ which, when all pressed, illuminates those bulbs.  It follows that pressing some combination of buttons on $\psi$ illuminates the configuration $\II(g)$.
\end{proof}

\ms
\begin{lemma} \label{nine}
All $g \in K$ with $\set{-n,n}\subseteq\II(g)$ and $\LL(g)=\0$
satisfy $d_{\AA}(g,1) \geq 4n$.
\end{lemma}

\begin{proof} Retracting $K$ onto $G$ by identifying $t$ and $u$
sends $\AA$ to $\SS= \set{a,s,t,at,ta,ata,as,sa,asa}$ (with some repetitions) and sends $g$ to an element $\bar{g} \in G$ satisfying the hypotheses of Proposition~4 of  \cite{CR}.  So $d_{\SS}(1, \bar{g}) \geq 4n$, by Proposition~4 of  \cite{CR}, and it follows that $d_{\AA}(g,1) \geq 4n$.
\end{proof}

\ms

Lemmas~\ref{distmax} and \ref{nine} imply $(K,\AA)$ has unbounded
 depth: all $g_n\in K$ with $\II(g_n)=\set{-n,n}$ and
 $\LL(g_n)=\0$ satisfy $d_{\AA}(1,g) =4n$ by the lemmas, and for all group elements
$h$ with $d_{\AA}(1,h) > 4n$ we have  $d_{\AA}(g_n,h) > n$  steps since
$\LL(h) \notin P_n$ by Lemma~\ref{distmax} and it takes at least $n$ steps to get outside $P_n$ from $\mathbf{0}$. So the depth of $g_n$ in $(K,\AA)$ must be at least $n$. 

\ms

We now turn to the other half of Theorem~\ref{Thm K}.  
The following definition is similar to that of a $\DD$-path in Section~\ref{proof of Thm H}.

\begin{definition}
We define a \emph{$\BB$-path} $\mu$ in $\Reals^3$ to be a path that runs from $\mathbf{0}$ to some point in $\Z^3$ and is obtained by concatenating two kinds of line-segments:
\begin{itemize}
\item slope-$\pm1$
\textup{(}\emph{diagonal}\textup{)} segments of length $1/\sqrt{2}$ parallel to the $q,r$-plane
and connecting points in $\Z^3$ to points in $\Z\times(1/2+\Z)^2$, and
\item
\textup{(}\emph{vertical}\textup{)} segments of length $1$ perpendicular to that plane
and connecting either points of $\Z^3$ to each other \textup{(}that is lying at integer points of the $q,r$-plane\textup{)} or points of
$\Z\times(1/2+\Z)^2$ to each other \textup{(}that is, lying at half-integer points of that plane\textup{)}.
\end{itemize}
\end{definition}

Unfortunately, the relationship between $\BB$-paths and words on $\BB^{\pm 1}$ is not as straightforward as that between $\DD$-paths and words on $\DD^{\pm 1}$, since $\BB$-paths can include vertical segments where the $q'$ and $r'$ co-ordinates are half-integers.  To begin to handle the subtleties we make the following definitions.

\begin{definition}
A \emph{decorated $\BB$-path} is a pair $(\mu, \mathcal{V})$ consisting of a  $\BB$-path $\mu$ and a set $\mathcal{V}$ of vertices along $\mu$ \textup{(}at endpoints of segments\textup{)}. If $\mu$ runs from $\mathbf{0}$ to $\LL(h)$ and pressing all the buttons at vertices in $\mathcal{V}$ en route illuminates the bulbs $\II(h)$ then we say $(\mu, \mathcal{V})$ \emph{represents} $h$.  
\end{definition}

\begin{definition}
A subpath of a \textup{(}possibly decorated\textup{)} $\BB$-path is a \emph{$\Z^3$-subpath} if it meets $\Z^3$ only at its end points.   \textup{(}So a \emph{$\Z^3$-subpath} in a $\BB$-path is either a vertical segment between two integer points or a pair of diagonal segments separated by some number of vertical segments, each connecting  half-integer points.\textup{)}
\end{definition}

\begin{definition}
A \emph{decorated} $\BB$-path $(\mu, \mathcal{V})$ is \emph{word-like} if the projections of all the vertical segments in $\mu$ to the $q,r$-plane are to points in $\Z^2$.   
\end{definition}

\ms

Translation between words $w$ on $\BB^{\pm 1}$ and word-like decorated $\BB$-paths  $(\mu,\mathcal{V})$, representing the same element of $K$, works as follows.  The key is that if $(\mu,\mathcal{V})$ is word-like then each of the $\Z^3$-subpaths in $\mu$ is either a vertical segment or is a concatenation of two diagonal segments.

\ms

\ni \emph{Obtaining $w$ from $(\mu,\mathcal{V})$.}    

\ss
From a word-like decorated $\BB$-path we read $w$ off in a  similar way to how we obtained words from $\DD$-paths.    Each $\Z^3$-subpath in $\mu$ corresponds to a generator or an inverse generator: vertical segments to $s^{\pm1}$, and pairs of diagonal segments to $t$, $u$, $tu$, $ut^{-1}$, $t^{-1}t$, $tt^{-1}$ or one of their inverses, in each case interspersed with letters $a$ as required to press the buttons at vertices in $\mathcal{V}$.  Let $w$ be the word on $\BB^{\pm 1}$ which, when read from right to left, has letters corresponding to the $\Z^3$-subpaths of $\mu$ (beginning at the start of $\mu$).   

\ms

\ni \emph{Obtaining  $(\mu,\mathcal{V})$ from  $w$.}    

\ss
Read  $w$ from right-to-left to obtain $(\mu, \mathcal{V})$ as follows.  Ignoring all $a$, let the $s^{\pm1}$ give vertical segments and the $t$, $u$, $tu$, $ut^{-1}$, $t^{-1}t$, $tt^{-1}$ (and their inverses) give appropriate pairs of diagonal segments, advancing $\mu$ in a manner corresponding to their actions; the locations of the $\mathcal{V}'$ are then dictated by the positions of the $a$'s in $w$.  

\ms

One might hope that, in analogy with the proof of Theorem~\ref{Thm H}, given $g \in K$, there exists a word-like decorated $\BB$-path $(\mu, \mathcal{V})$,  from which the word obtained is a geodesic on $\BB^{\pm1}$ representing $g$, and all the vertices $\mathcal{V}$ are on an initial segment of  $\mu$ that does not leave the plane $q=r$.  The truth, as presented in the following lemma, can be marginally more complicated (but only when the final vertex of $\mu$ in the plane $q=r$ has half-integer $q$- and $r$-co-ordinates).  

\begin{lemma}\label{one and a half}
Given $g \in K$, there exists a word-like decorated $\BB$-path $(\mu, \mathcal{V})$ with the following properties.  Firstly, the word obtained from $(\mu, \mathcal{V})$ is a geodesic on $\BB^{\pm1}$ representing $g$.  Secondly, the $p$-co-ordinate either monotonically increases or monotonically decreases along $\mu$.  Thirdly, $\mu$ is a concatenation of four arcs $\mu_1, \mu_2, \mu_3, \mu_4$, such that $\mu_1$ is in the $q=r$ plane,  $\mu_2$ is at most one diagonal segment perpendicular to the $q=r$ plane, $\mu_3$ is a concatenation of vertical segments with no backtracking, and $\mu_4$ is a concatenation of diagonal segment perpendicular to the $q=r$ plane with no backtracking.  Fourthly, $\mu_2$, $\mu_3$ and $\mu_4$ are all on the same side of the $q=r$ plane.  And finally, all the vertices $\mathcal{V}$ are on $\mu_1$, $\mu_2$ and $\mu_3$.   
\end{lemma}

\begin{proof}
Take any geodesic word $w_0$ on $\BB^{\pm 1}$ representing $g$ and obtain from it a word-like decorated $\BB$-pair $(\mu, \mathcal{V})$.  Let $D$ and $V$ be the number of diagonal and vertical segments in $\mu$, respectively.  We will alter $(\mu, \mathcal{V})$ in three steps until it satisfies the conditions of the lemma.  At all times during the transformation it will continue to represent $g$. Furthermore, neither $V$ nor $D+V$ will increase. Since a word-like path with $V$ vertical segments and $D$ diagonal segments corresponds to a word of length $V+D/2$, the resulting word will still be geodesic.

\ss

\ni (\emph{i}) \emph{Collect all diagonal segments $\sigma$ in $\mu$ running perpendicular to the plane $q=r$ at the end of $\mu$.} 

\ss
Removing such a $\sigma$, translating the portion of $\mu$ after $\sigma$ and the vertices of $\mathcal{V}$ thereon parallel to $\sigma$ to close up the gap, and then re-attaching $\sigma$, produces a new decorated $\BB$-path representing $g$ (as all vertices in $\mathcal{V}$ that are moved are replaced with vertices at the same values of $p$ and $q+r$).  

Repeat until all such $\sigma$ are collected at the end. So some initial segment of $\mu$ runs in the plane $q=r$ and then the remaining terminal segment of $\mu$ runs perpendicular to the plane $q=r$.   Move  the vertices $\mathcal{V}$ so that they lie on this initial segment, and remove any backtracking from the terminal segment.  

Note that $V$ has not changed and $D$ has not increased.  Moreover, $\mu$ now lies entirely in one of the two closed half-spaces bounded by the plane $q=r$.

\ss

It may be that $(\mu, \mathcal{V})$ is no longer word-like.  We rectify this with our next two steps.

\ss

\ni (\emph{ii}) \emph{Ensure the $p$-co-ordinate changes monotonically along $\mu$.}   

Suppose $\nu$ is a subpath of $\mu$ that begins with a vertical segment along which $p$ increases by $1$, continues with a number of diagonal segments, and concludes with a vertical  segment along which $p$ decreases by $1$.  Furthermore suppose the $p$-coordinate is non-negative at every point of along $\nu$.  (Note that $\nu$ must be in the plane $q=r$.)  Let $\hat{p}$ be the $p$-co-ordinate of the initial point of $\nu$, and let $m$ and $M$ be the minimum and maximum values of $q$ at points $(p,q,q)$ on $\nu$.   Remove $\nu$ from $\mu$ and in its place insert a path $\hat{n}$ that runs along the line $p=\hat{p}$ of the $q=r$ plane and visits $(\hat{p},m,m)$ and  $(\hat{p},M+1/2,M+1/2)$.  Such a $\hat{\nu}$ exists with no more than two more diagonal segments than $\nu$ and (obviously) two fewer vertical segments.  Moreover, vertices of $\mathcal{V}$ on $\nu$ can be replaced by vertices on $\hat{\nu}$ in such a way that the lights illuminated does not change -- we leave the details of this to the reader save to say that the reason $\hat{\nu}$ is made to overshoot $\nu$ in the $(0,1,1)$-direction is so that the effect of pressing a button at a point $(\hat{p}+1,q,q)$ of maximal $q$ on $\nu$ can be duplicated one layer lower.   

A similar result holds for subpaths along which the $p$-coordinate is non-positive.  

It follows that, after exhaustively making such changes, $(\mu, \mathcal{V})$ can be made to have the $p$-co-ordinate either increasing or decreasing monotonically along $\mu$.  Moreover, in the course of step (\emph{ii}), any increase in $D$ is compensated for by at least the same decrease in $V$, so the total number of segments does not increase.

\ss

\ni (\emph{iii}) \emph{Ensure there is no $\Z^3$-subpath $\xi$ in $\mu$ that is made up of a diagonal segment $\xi_0$ followed by a concatenation of vertical segments $\xi_1$, all above a half-integer point in the $q,r$-plane, and then another diagonal segment $\xi_2$.}

\ms 

Suppose there is such $\Z^3$-subpath $\xi$.  We explain how to change $\mu$ to remove it.  We assume that $\mu$ lies on the nonnegative-$p$ side of the $q,r$-plane.  The case where $\mu$ lies on the nonpositive-$p$ side is similar, and we omit it.  

Suppose first that $\xi_0$ and $\xi_2$ are in the plane $q=r$. So motion along them increments or decrements $q+r$.     Note that by (\emph{ii}), the $p$-co-ordinate is monotonically increasing along $\xi_1$.  

\begin{figure}[ht]
\psfrag{p}{$p$}
\psfrag{q+r}{$q+r$}
\centerline{\epsfig{file=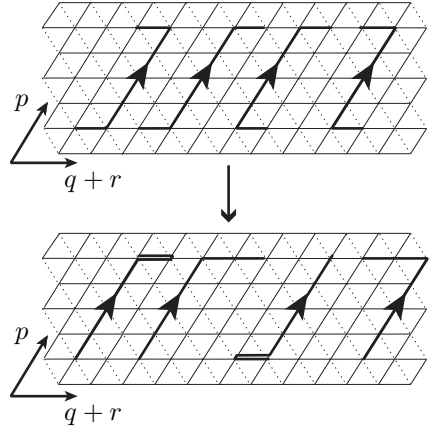}} \caption{Pushing the subpath $\xi$ of $\mu$ -- four cases} \label{pushing}
\end{figure}

Regard the plane $q=r$ as the rhombic grid of the lamplighter model for $G$ as in Figure~\ref{pushing}.  If $q+r$ increases along $\xi_0$ (the first and second cases shown in the figure) then interchange $\xi_0$ and $\xi_1$ -- that is, push $\xi$ across a parallelogram, replacing its lower side  ($\xi_0$) and its right side ($\xi_1$) by its left and upper sides.   Failing that, if $q+r$ increases along $\xi_2$ (the third case in the figure) then interchange $\xi_1$ and $\xi_2$ -- that is, push across the parallelogram with left side ($\xi_1$) and top side ($\xi_2$).  The remaining possibility (the fourth case in the figure) is that $q+r$ decreases along both $\xi_0$ and $\xi_2$.  In this case we exchange $x_0$ and $x_1$ -- in other words we push $\xi$ across the parallelogram with $x_0$ as its lower side and $x_1$ as its left side, replacing those two sides with the right and top.  We leave it for the reader to check that in every case $\mathcal{V}$ can be altered, selecting new vertices on the other side of the parallelogram, illuminating the same bulbs -- perhaps the easiest way to do this is push $\xi$ across the parallelogram one unit-sided rhombus at a time, working from bottom to top, checking that at each stage vertices can be moved across in a way that maintains the same configuration of illuminated lights.  

Suppose next that $\xi_0$ and $\xi_2$ are not both in the plane $q=r$. The only way this can happen is for $\xi_2$ to be the first segment of the portion of $\mu$ outside the plane $q=r$.  In this case interchange $\xi_1$ and $\xi_2$, translating any vertices of $\mathcal{V}$ on $\xi_1$ on the way. Since motion perpendicular to the plane $q=r$ does not affect which lights are toggled by the buttons, the new path will represent the same group element. 

Repeating this process eventually renders $(\mu, \mathcal{V})$ word-like and it changes neither $D$ nor $V$.  If the interchange in the paragraph above was implemented then the vertical path moved is $\mu_3$, the diagonal segment that now precedes it is $\mu_2$, and the remaining initial and terminal portions of $\mu$ are $\mu_1$ and $\mu_4$.  Otherwise, $\mu_2$ and $\mu_3$ are taken to be of zero-length and $\mu_1$ is the maximal length initial path of $\mu$ lying in  the plane $q=r$ and $\mu_4$ is the remainder.  

\ms

Let $w$ be the word now obtained from $(\mu, \mathcal{V})$.  The length of $w$ is $V+ D/2$, which is at most the length of $w'$.  So $w$ is a geodesic.  Moreover, all vertices of $\mathcal{V}$ are on $\mu_1$, $\mu_2$ and $\mu_3$, as required.
\end{proof}

We can now conclude our proof. Suppose $g \in K \ssm \set{a}$. Consider $g' \in K$ with $\II(g')=\II(g)$ and $\LL(g')$ located so that the line through $\LL(g)$ and $\LL(g')$ is perpendicular to the plane
$q=r$ and $\LL(g')$ is  a distance $\sqrt{2}$ further from that plane than is $\LL(g)$.  Then $g'$ is at distance greater than $1/\sqrt{2}$ from the plane $q=r$ and so, by Lemma~\ref{one and a half}, some geodesic word $w'$ for $g'$ comes from a decorated $\BB$-path $(\mu',\mathcal{V}')$ whose last two segments form the diagonal from $\LL(g)$ to $\LL(g')$, and vertices of $\mathcal{V}'$ lie beyond $\LL(g)$. So deleting these last two segments gives a decorated $\BB$-path $(\mu, \mathcal{V}')$ yielding a word $w$ representing $g$, which is a subword of $w'$ by construction. Since $d_{\BB}(g,g')=1$, our proof is complete.

\section{Dead-end depth in $\Z$} \label{depth in Z}
In this section we prove Proposition~\ref{depth in Z result}.

 Writing $r \in \set{0, 1, \ldots,n(n+1)-1}$ as $an+b$ for some $0\leq
 a \leq n$ and $0\leq b<n$, we have $$r \ = \ (a-b)n+b(n+1) \ = \
 (a-b+n+1)n + (b - n)(n+1).$$ Thus if $d$ denotes the word metric on
 $(\Z,\set{n,n+1})$ then $$d(0,r) \ \leq \ \min \set{ \abs{a-b} + \abs{b}, \abs{ a -b+n+1 } + \abs{b-n} },$$ which is at most $n$ because if $a+n \geq 2b$ then $\abs{a - b} + \abs{b} \leq n$ and if $a +n < 2b$ then $$ \abs{a -b +n +1} + \abs{b-n} \ = \ (a-b+n+1) + (n-b) \ = \ a -2b +2n +1 \ < \  n+1,$$ and so $\abs{a -b +n +1} + \abs{b-n}  \leq \ n$.
 
Thus $$\max \set{ \ d(0,r) \  \mid \  r =0,1, \ldots, n(n+1) -1 \ } \ \leq \ n. $$  It is enough to show that this maximum is attained at a positive integer $<n+1$, because then any point of
$\Z$ further from $0$ will be at least $n(n+1)-n=n^2$ away in the
standard metric on $\Z$, and hence at
 least $n^2/(n+1)>n-1$ away in $d$.  
 
 If $n$ is even, we set $n=2m$;
the maximum will be attained at $m$ since if $m=a(2m)+b(2m+1)$ then $$m \
= \ b \ \ \ (\textup{mod} \ 2m), \ \ \ \ \ \ \ \ m \ = \ -a \ \ \
 (\textup{mod} \ 2m+1)$$ and so $\abs{a}$, $\abs{b} \geq m$. If $n$ is
odd, we set $n=2m-1$;  the maximum
 will again be attained at $m$ since if $m=a(2m)+b(2m-1)$ we have $$m \ =
\ -b \ \ \ (\textup{mod} \ 2m), \ \ \ \ \ \ \ \ m \ = \ a \ \ \
 (\textup{mod} \ 2m-1),$$ and so $\abs{b} \geq m$ and $\abs{a} \geq m-1$.
In both cases $\abs{a} + \abs{b} \geq n$ and $m < n+1$.  This completes the proof.

\ms
We remark that Proposition~\ref{depth in Z result} also holds for groups $G$ such that $G/N \cong \Z$ for some finite $N \trianglelefteq G$ and for groups $\Z \rtimes H$ where $H$ is finitely generated. Whether it holds for other groups (even a rank $2$ free group) remains open.

\bibliographystyle{plain} \bibliography{bibli}

\ni \textsc{T.R.Riley, Mathematics Department, 310 Malott Hall, \\ Cornell University, Ithaca, NY 14853-4201, USA}  \\ \texttt{tim.riley@math.cornell.edu} 

\ms

\ni \textsc{A.D.Warshall, Mathematics Department, 10 Hillhouse Avenue, \\ P.O. Box
208283, New Haven, CT 06520-8283, USA} \\ \texttt{andrew.warshall@yale.edu}

\end{document}